\documentclass[aps, prd, onecolumn, amssymb, amsmath, amsfonts]{revtex4-2}
\usepackage{amsthm}
\usepackage{graphicx}
\usepackage{hyperref}[hidelinks]

\def\fracs#1#2{{}^{#1}\hspace{-0.25em}/%
	\hspace{-0.15em}_{#2}^{\vphantom{1}}}

\makeatletter
\gdef\th@remark{%
	\thm@headfont{\small\scshape}%
}
\makeatother
\theoremstyle{remark}
\newtheorem{rem}{Remark}

\DeclareMathOperator\arctanh{arctanh}

\begin{document}
\title{Classification of the trajectories  of uncharged particles
in the Schwarzschild-Melvin metric}
\author{Ivan Bizyaev}
\affiliation{Ural Mathematical Center, Udmurt State University, Universitetskaya 1, 426034 Izhevsk, Russia}
\email{ bizyaevtheory@gmail.com}

\begin{abstract}
This paper investigates the trajectories of neutral particles in the
Schwarzschild-Melvin spacetime.
	  After reduction by cyclic coordinates this problem reduces to investigating a
two-degree-of-freedom Hamiltonian system that has no additional integral. A classification of regions of possible motion
of a particle is performed according to the values of the momentum and energy integrals.
Bifurcations of periodic solutions of the reduced system are analyzed using a Poincar\'{e}
map.

\end{abstract}
	
	\maketitle

\tableofcontents
\section{Introduction}
	A qualitative analysis of the influence of a magnetic field on physical processes in a
neighborhood of a black hole is made possible by exact solutions to the Einstein-Maxwell
equations.
	In this paper we consider a solution, obtained by Ernst \cite{Ernst1976},
	which describes a nonrotating black hole immersed into an external magnetic field.
It generalizes the Schwarzschild solution and is a unique axisymmetric static solution 
(see \cite{HiscockBlack1981} for details) which, at a
large distance from the event horizon, becomes the magnetic Melvin universe
\cite{Melvin1964}.  In \cite{Wild1980Surface}, the Gaussian curvature of the event horizon
is calculated for this Schwarzschild-Melvin solution,
	 and its linear stability is considered in \cite{Brito2014Superradiant}.
A detailed description of the properties of the Schwarzschild-Melvin spacetime is
given in \cite{StuchlikPhoton1999}.
	
	Recently, zero geodesics for the Schwarzschild-Melvin metric have been considered in \cite{Junior2021Shadows}
  to describe the shadow of the black hole and the gravitational
lensing. The motion of charged particles was examined in \cite{Lim2015Motion}.
We note that, in investigating the motion of charged particles, one can sometimes
neglect the influence of the magnetic field on the metric
(see, e.g., \cite{Frolov2010Motion}). In this case it is necessary to find
a solution to the Maxwell equation for the fixed metric. For example,
an analog of the homogeneous magnetic field directed along the symmetry axis for the
Kerr metric was found in \cite{Wald1974Black}.

In this paper, a detailed analysis is made of the trajectories of neutral particles in the Schwarzschild-Melvin spacetime, which are described by time-like geodesics. The aim of this paper is not to discuss specific astrophysical objects with a strong magnetic field, but rather to study the dynamics of geodesics using, among others, various methods of qualitative analysis from dynamical systems theory.

	The equations of the geodesics of the Schwarzschild-Melvin metric possess a discrete symmetry associated with reflection relative to the equatorial plane. Their fixed points define the invariant submanifold on which plane trajectories lie. These trajectories are explored in detail in \cite{Galtsov1978Black}.  
There it is shown that circular trajectories of light beams can only exist if the magnetic field does not exceed a certain value.
In Section \ref{HillReg} it is shown that this result persists for neutral particles,  earlier this result was obtained in \cite{Esteban1984}. 
The motion of charged particles in a neighborhood of a weakly magnetized Schwarzschild black hole was considered in \cite{Frolov2010Motion}. There an analysis of the trajectories lying in the equatorial plane is made and it is shown that the amplification of the magnetic field leads to the emergence of a curly cycloid-like trajectory. 	

Circular trajectories define relative equilibria in the initial system because they are fixed 
points of reduced equations of motion. In their turn, the relative equilibria coincide with the critical points of the Hamiltonian, which in this case define the critical points of the effective potential. As a 
result, the analysis of the relative equilibria allows one to classify regions of possible 
motion according to the values of first integrals (see \cite{BolsinovBifurcation2012, Smale1970Topology, SmalePlanar1970} for details).
This classification, presented in Section  \ref{HillReg}, singles out the values of
angular momentum and energy integrals, as well as those of the magnetic field in which 
all trajectories are bounded.

In Section \ref{regularization}, a regularization of the reduced equations of motion is 
performed. Next, in Section  \ref{sectionMap}, the bifurcations of the periodic solutions of the regularized reduced system are considered using a Poincar\'{e} map. In particular, it is shown that, for a sufficiently strong magnetic field in which there are no circular trajectories, other bounded trajectories can occur. Also, in Section \ref{sectionMap}, an example of such a bounded trajectory arising after a saddle-node bifurcation is given.

 In contrast to the Schwarzschild-Melvin metrics, the equations of motion
 for the geodesics under consideration are in the general case nonintegrable by quadratures.
	This conclusion is drawn in \cite{Li2019Chaotic, Wang2017Chaos, Karas1992Chaotic} using a numerically constructed Poincar\'{e} map which shows chaotic trajectories. Then, the largest Lyapunov exponent is calculated for these trajectories and it is shown that it takes a positive value. It should be noted that such a numerical analysis  is made for specific parameter values and does not lead to the conclusion about nonintegrability of the system for various parameter values.
One of the efficient tools for analytically proving the absence of an additional integral in the system is to calculate the Melnikov integral since it allows one to show the intersection of separatrices \cite{KozlovIntegrability1983, RobinsonDynamical1995}. In Section \ref{separatrix}, for the system under consideration, the Melnikov integral is calculated and the intersection of separatrices is shown. Next, a visualization of these separatrices in the Poincar\'{e} map is presented using the numerical algorithm proposed in \cite{Li2012Algorithm}.

\section{The Schwarzschild\,--\,Melvin metric}
The Schwarzschild-Melvin spacetime is a static and axisymmetric solution
to the Einstein-Maxwell equations. In the Schwarzschild coordinates $\boldsymbol{x}=(t, r, \theta, \varphi)$, its line element can be represented in the following form \cite{Ernst1976}:
\begin{equation}
\begin{gathered} 
\label{eq_g}
ds^2 = g_{ij}dx^idx^j = \Lambda^2\left[ - \left(1 - \frac{r_s}{r} \right)c^2dt^2 +  \frac{dr^2}{1 - \dfrac{r_s}{r} }+ r^2d\theta^2  \right]  +\frac{r^2}{\Lambda^2}\sin^2 \theta d\varphi^2,\
\Lambda = 1 + \frac{B^2}{4}r^2\sin^2\theta,
\end{gathered}
\end{equation}	
where summation is over repeating indices, $c$ is the velocity of light, and the signature $(-,+,+,+)$ is chosen.

 The electromagnetic potential is the following 1-form:
\begin{equation}
\label{eq_Ai}
A_idx^i = - \frac{c^2Br^2}{2\sqrt{G}\Lambda}\sin^2\theta d\varphi,
\end{equation}
where $G$ is the gravitational constant.
Relations \eqref{eq_g} and \eqref{eq_Ai} contain two parameters $r_s$ and $B$. For them, one can single out two special cases.
\begin{itemize}
\item[1.] {\it Schwarzschild's solution} \cite{Schwarzschild1916} ($B=0$, $r_s\neq 0$), which describes an isolated spherically symmetric black hole with the event horizon
$$
\mathcal{S}_h = \{ (t, r, \theta, \varphi) \ | \ r=r_s \},
\quad r_s = \frac{2MG}{c^2},
$$
where $M$ is the mass of the black hole.
\item[2.] {\it Melvin's solution} \cite{Melvin1964} ($B\neq0$, $r_s =0$), which corresponds to the space filled by a static magnetic field. Its vector of magnetic induction always has the same direction, but takes different absolute values. As $r\to +\infty$, the magnetic field does not tend to zero, and the
    metric \eqref{eq_g} is not asymptotically flat.
\end{itemize}
If $r_s\neq 0$ the metric \eqref{eq_g} always possesses the event horizon $\mathcal{S}_h$.
Therefore, the solution under consideration can be interpreted as a black hole immersed
into the external magnetic field. The nonzero components of the magnetic induction of this
field have the form \cite{Ernst1976}
\begin{equation}
\label{eq_B}
B_r=\frac{c^2B}{\sqrt{G}\Lambda^2}\cos\theta, \quad B_\theta = - \frac{c^2B}{\sqrt{G}\Lambda^2}\sqrt{ 1 - \frac{r_s}{r}}\sin\theta.
\end{equation}

Let us make a qualitative estimate of what magnitude the magnetic field \eqref{eq_B} can have.
For this, we note that the parameter $B$ has an inverse length. Therefore, it is convenient
to introduce a dimensionless parameter $B_s=B r_s$. From this we obtain the following
typical scale of the strength of the magnetic field  \eqref{eq_B}:
$$
B_M \sim \frac{c^4 B_s}{G^{3/2}M} \sim 10^{19} \frac{B_s M_\odot}{M} \ \mbox{Gauss},
$$
where $M_\odot$ is the mass of the Sun.
In a neighborhood of a black hole the magnetic field is assumed to be caused by the accretion disk of the black hole. Specifically, for supermassive black holes with a mass of order $M \sim 10^9 M_\odot$
a magnetic field $B_M \sim 10^4$ Gauss \cite{Piotrovich2010Magnetic} is observed, which gives $B_s \sim 10^{-6}$.
This magnetic field does not influence greatly the geometry of spacetime in a neighborhood
of the event horizon  $\mathcal{S}_h$, but can do so when $r \gg r_s$.

A magnetic field capable of influencing the geometry of spacetime in a neighborhood of the event horizon can arise  after merging of the black hole with a magnetar (for details, see \cite{Chakraborty2024Black} and references therein). A magnetar is a neutron star which possesses a strong magnetic field. For example, the magnitar SGR J1745-29 with a magnetic field of order $10^{14}$ Gauss is situated near Sagittarius A* (Sgr A*) at the center of the Milky Way.

We note that the solution \eqref{eq_g}-\eqref{eq_Ai} does not take into account a number of important factors, for example, the rotation of the black hole. Nevertheless, it is of interest what influence the magnetic field described by the Schwarzschild-Melvin solution can have on the particles' trajectories.
Throughout the remainder of this paper, the analysis of these trajectories will be performed using a system of units in which the gravitational constant and the velocity of light are equal to 1 ($G=1$, $c=1$). 
In addition, in the metric \eqref{eq_g} we transform to the dimensionless variables and  the dimensionless parameter
$$
r \to r_s r, \quad  t \to r_s t, \quad B \to \frac{B}{r_s},
$$
then the event horizon corresponds to the value $r=1$. Next, in constructing the figures for the dimensionless parameter we choose
 $B\in[0.01, 0.4]$, which corresponds to a fairly strong magnetic field.

 The electromagnetic tensor $F_{ij}$ contains nonvanishing components \cite{Biak1985}:
$$
F_{\theta \varphi} = \frac{B}{\Lambda^2} r^2\sin \theta\cos\theta, \quad F_{r \varphi}=\frac{B}{\Lambda^2}\sin^2\theta.
$$

\begin{figure*}
\includegraphics[scale=0.8]{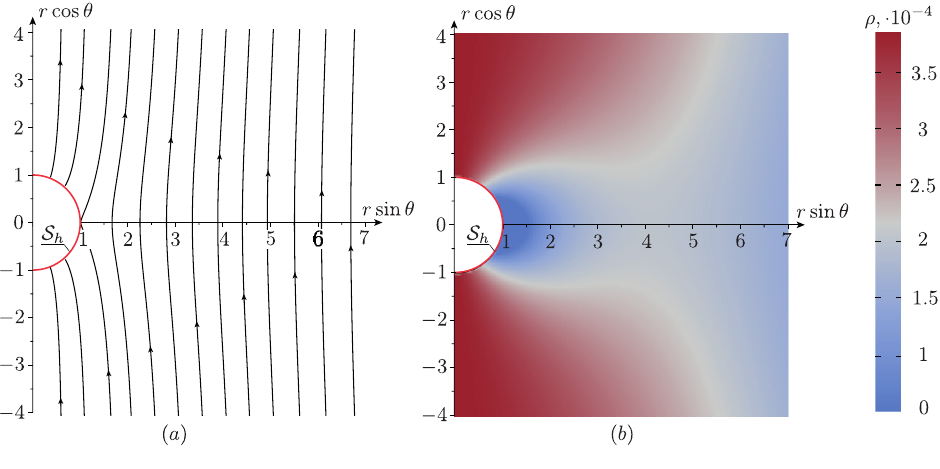}
\caption{On the plane ($r\sin \theta, r \cos \theta$) for fixed $r_s=1$ and $B=0.1$ the figure shows: (a) force lines of the magnetic field; (b) the value of the energy density of the magnetic field~$\rho$. }
\label{fig01}
\end{figure*}

The energy density of the magnetic field \eqref{eq_B} can be determined as follows:
$$
\rho=-\frac{1}{16\pi} F_{ij} F^{ij}=\frac{B^2}{16\pi\Lambda^4}\left( 1 - \frac{\sin^2\theta}{r} \right)
$$
Figure \ref{fig01} shows the force lines of the magnetic field and the value of density $\rho$. As can be seen,
the magnetic field takes the largest value in a neighborhood of the symmetry axis and the smallest value in the equatorial plane in a neighborhood of the horizon $\mathcal{S}_h$.

\begin{rem}
In the general case, the Schwarzschild-Melvin solution contains another parameter, $\psi$, on which only the vector potential \cite{Griffiths2009} depends. In the chosen dimensionless variables it has the following form:
$$
A_idx^i = B\cos\psi (r-1)\cos\theta dt - \frac{Br^2}{2\Lambda}\sin \psi \sin\theta^2d\varphi.
$$
If $\psi= 0$, then the previous potential describes the electric field, and for $\psi=\dfrac{\pi}{2}$
we obtain the magnetic field \eqref{eq_Ai}.
\end{rem}

Consider the motion of a material point (particle) only on the ``external'' side of the
event horizon $\mathcal{S}_h$, i.e., on the manifold
$$
\begin{gathered}
\mathcal{N}^4=\{ (t, r, \theta, \varphi) \ | \ t\in(-\infty, +\infty), \ r\in(1, +\infty), \\ \theta\in(0,\pi), \ \varphi\in(0, 2\pi)  \}.
\end{gathered}
$$
Let the electric charge of the particle be zero. Then its trajectory is a geodesic on
$\mathcal{N}^4$.
We use a natural parametrization of the trajectory of particle $\boldsymbol{x}(\tau)$, where $\tau$ is the proper time of the particle.
Then the required trajectories satisfy equation \cite{Misner1973} ---
\begin{equation}
\label{eqH}
\begin{gathered}
\frac{d\boldsymbol{x}}{d\tau} = \frac{\partial H}{\partial \boldsymbol{p}}, \quad \frac{d\boldsymbol{p}}{d\tau} = - \frac{\partial H}{\partial \boldsymbol{x}}, \\
H=\frac{1}{2}g^{ij}p_ip_j,
\end{gathered}
\end{equation}
where $\boldsymbol{p}=(p_t, p_r, p_\theta, p_\varphi)$ is the momentum of the particle and $g^{ij}$ is the matrix inverse to the metric defined by relation \eqref{eq_g}.

The metric $g_{ij}$ does not depend explicitly on time $t$ and angle $\varphi$,  therefore, they are cyclic coordinates for equations \eqref{eqH}. As a consequence, the corresponding momenta remain unchanged:
$$
E=-p_t={\rm const}, \quad L = p_\varphi = {\rm const}.
$$
From a physical point of view $E>0$ is the energy of the material point and $L$ is the projection of its angular momentum onto the symmetry axis.

Thus,  the Hamiltonian system with two degrees of freedom decouples from the system \eqref{eqH}:
\begin{equation}
\label{eqRed}
\begin{gathered}
\frac{dp_r}{d\tau} = - \frac{\partial H}{\partial r}, \quad \frac{dp_\theta}{d\tau} = - \frac{\partial H}{\partial \theta}, \quad
\frac{dr}{d\tau} = \frac{\partial H}{\partial p_r}, \quad \frac{d\theta}{d\tau} = \frac{\partial H}{\partial p_\theta}, \\
H= \frac{\Gamma}{2\Lambda^2}p_r^2 + \frac{p_\theta^2}{2r^2\Lambda^2} + U(r, \theta), \quad
U(r,\theta) = \frac{L^2\Lambda^2}{2r^2\sin^2\theta} -   \frac{E^2}{2\Lambda^2\Gamma}, \quad   \Gamma=1 - \frac{1}{r},
\end{gathered}
\end{equation}
and the phase space of this system has the form
$$
\mathcal{M}^4 = \{ \boldsymbol{y}=(p_r, p_\theta, r, \theta) \ | \ r\in(1, +\infty), \ \theta\in(0, \pi) \}.
$$
The trajectories of the material particles lie on the fixed level set of the Hamiltonian \cite{Misner1973}
\begin{equation}
\label{eq_h_level}
H(\boldsymbol{y})=-\dfrac{1}{2}.
\end{equation}

According to the known solution of the system \eqref{eqRed} the evolution of the remaining variables is defined from the equations
\begin{equation}
\label{eqFull}
\frac{d \varphi}{d\tau} = \frac{\Lambda^2 L}{r^2 \sin^2 \theta}, \quad
\frac{dt}{d\tau} = \frac{E}{\Lambda^2\Gamma}.
\end{equation}

In the case $L=0$ the system \eqref{eqRed} contains no singularity at the poles $\theta=0$ and $\theta=\pi$. As a result, the trajectories can reach these poles or lie entirely on the symmetry axis (see, e.g., \cite{StoghianidisPolar1987}).  This case requires a separate analysis; therefore, in what follows we will assume $L\neq 0$.

We note that the system of equations \eqref{eqRed} and \eqref{eqFull} is invariant under the transformation
$$
L\to -L, \quad \varphi \to -\varphi,
$$
therefore, in what follows, it can be assumed without loss of generality that $L>0$.
Moreover, the parameter $B$ appears in the system of equations \eqref{eqRed} and \eqref{eqFull}
only in the second degree. Therefore, we will also assume that $B>0$.

\section{Relative equilibria and the Hill regions}
\label{HillReg}
From a mechanical point of view the reduced system \eqref{eqRed} is a natural Hamiltonian
system with two degrees of freedom in which the function $U(r, \theta)$ is an effective
potential. In this case, for the fixed points of this system the momentum components
vanish: $p_r=0$, $p_\theta=0$, and the  coordinates take the fixed values $r=r_c={\rm const}$, $\theta=\theta_c={\rm const}$ and are critical points of the function $U(r, \theta)$, i.e., they satisfy  the
system
\begin{align}
\label{eqR_Th}
\left.\frac{\partial U}{\partial r}\right|_{\substack{r=r_c \\ \theta=\theta_c}}&=\frac{E^2}{2\Lambda_c^2(r_c-1)}\left( 3 - \frac{4}{\Lambda_c} + \frac{r_c}{r_c - 1} \right)  - \frac{\Lambda_c(2 - \Lambda_c)}{r_c^3\sin^2\theta_c}L^2=0, \nonumber \\
 \left.\frac{\partial U}{\partial \theta}\right|_{\substack{r=r_c \\ \theta=\theta_c}}&= \frac{2E^2r_c(\Lambda_c - 1)}{\Lambda_c^3(r_c - 1)\sin\theta_c}\cos\theta_c - \frac{\Lambda_c(2 - \Lambda_c)}{r_c^2\sin^3\theta_c}L^2\cos\theta_c = 0, \\
&  \quad \Lambda_c = 1 + \frac{B^2}{4}r_c^2\sin^2\theta_c. \nonumber
\end{align}
For the solutions of this system the following equation is satisfied:
\begin{equation}
\label{eq_0}
r_c\cos \theta_c \left.\frac{\partial U}{\partial r}\right|_{\substack{r=r_c \\ \theta=\theta_c}} - \sin\theta_c  \left.\frac{\partial U}{\partial \theta}\right|_{\substack{r=r_c \\ \theta=\theta_c}} = \frac{r_cE^2\cos\theta_c}{2\Lambda_c^2(r_c - 1)^2}=0,
\end{equation}
which  has the unique solution
\begin{equation}
\label{eq_th0}
\theta_c=\dfrac{\pi}{2}.
\end{equation}
We note that equation \eqref{eq_0} can also be satisfied if $r_c\to +\infty$. This case
will be treated separately.
Thus, the fixed points of the system \eqref{eqRed} are {\it relative
equilibria} \cite{Marsden1974Reduction,Arathoon2023} at which the trajectory of the material point is a
circle lying in the equatorial plane. Consider them in more detail.

According to condition  \eqref{eq_h_level}, the critical points satisfy another
equation which, taking \eqref{eq_th0} into account, has the form
\begin{equation}
\label{eq001}
\begin{gathered}
U\left(r_c, \dfrac{\pi}{2}\right)= \frac{\left(1+\dfrac{B^2r_c}{4}\right)^2L^2 }{2r_c} 
- \frac{r_c E^2}{2(r_c - 1)}\left(1+\dfrac{B^2r_c}{4}\right)^{-2} =-\dfrac{1}{2}.
\end{gathered}
\end{equation}
We note that in the case   \eqref{eq_th0}  the second equation in the system \eqref{eqR_Th}
is satisfied identically, and the first equation can be represented as
\begin{equation}
\label{eq_05}
\begin{gathered}
\left(\frac{1}{r_c-1} + \frac{3}{r_c} - \frac{4}{br_c}  \right)U\left(r_c, \dfrac{\pi}{2}\right) + \frac{b L^2}{2r_c^3(r_c-1)} Z(r_c) =0, \\
Z(r_c) = 2r_c - 3 - B^2r_c^2\left( \frac{3}{2}r_c - \frac{5}{4}  \right), \\ b = 1 + \frac{B^2r_c^2}{4}.
\end{gathered}
\end{equation}

Equations \eqref{eq001} and \eqref{eq_05} depend in a fairly complicated way on the radial
coordinate $r_c$. At the same time, the values of the first integrals $L$ and $E$ appear in
\eqref{eq001} and \eqref{eq_05} in the second degree. Therefore, we consider these equations as a system
defining the values of the first integrals depending on $r_c$.
It can be shown that this system has a solution only for those values of $r_c$ that
satisfy the inequality $Z(r_c)>0$.
This inequality has a solution only for
\begin{equation}
\label{eq002}
0\leqslant B<B_{n}, \quad  B_{n}=\frac{4\sqrt{3}}{\sqrt{169+38\sqrt{19}}}\approx 0.379.
\end{equation}
If condition \eqref{eq002} is satisfied, then $r_c\in(r_1, r_2)$,  where $r_1$ and $r_2$
are the roots of the equation $Z(r_c)=0$:
\begin{equation}
\label{eq_r12}
\begin{gathered}
r_1 = \frac{1}{3}\sqrt{\frac{25}{9} + \frac{16}{B^2}}\cos \left(  \frac{\Delta}{3}  + \frac{4\pi}{3}\right) + \frac{5}{18}, \\
r_2 =  \frac{1}{3}\sqrt{\frac{25}{9} + \frac{16}{B^2}}\cos \frac{\Delta}{3}  + \frac{5}{18}, \\
 \Delta=\arccos \frac{B(125B^2 - 4752)}{(25B^2 + 144)^{\frac{3}{2}}}
\end{gathered}
\end{equation}
and
for $B=0$ we have $r_2 \to +\infty$.  Earlier, similar relations for  $r_1$, $r_2$ and $B_n$ --- but in a different system of units of measurement, were obtained in \cite{Esteban1984}.

\begin{figure*}[!ht]
\includegraphics[scale=0.85]{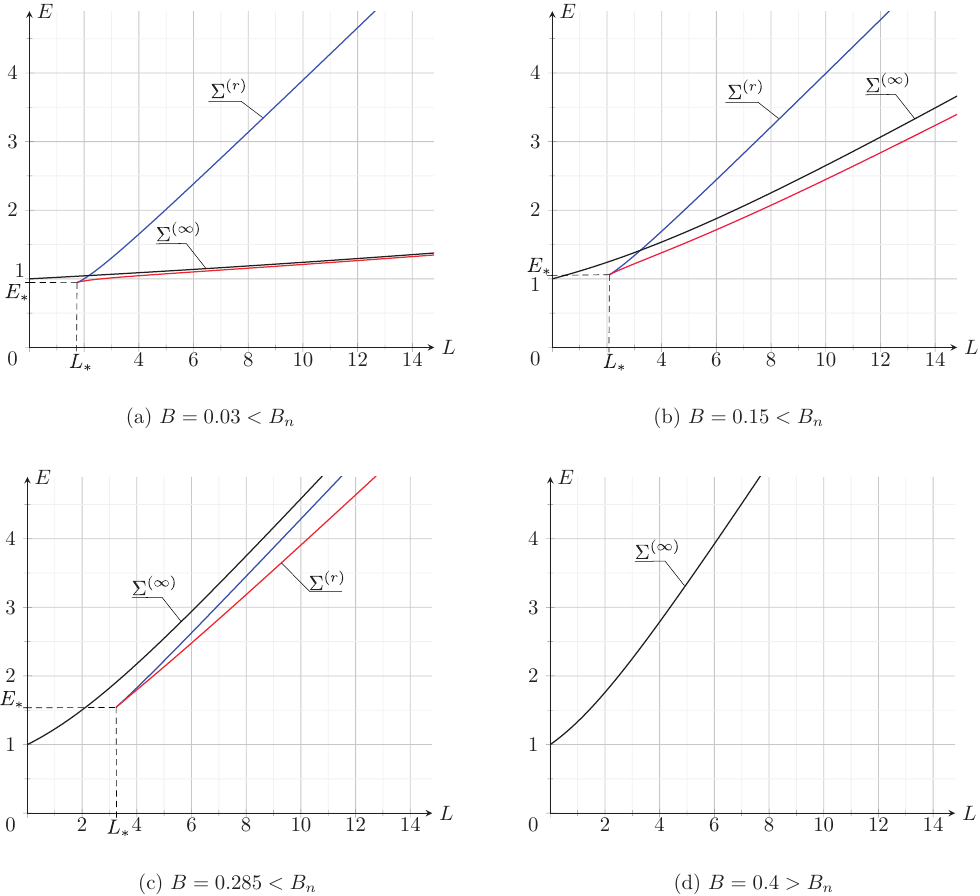}
\caption{Bifurcation diagram for different values of $B$. Red
indicates the part of the curve $\Sigma^{(r)}$ that corresponds to stable fixed points,
and blue indicates the part corresponding to unstable ones.
The point $(L_*, E_*)$ corresponds to a {\it cusp} type singularity on the
curve $\Sigma^{(r)}$, where $L_*=L_c(r_*)$ and $E_*=E_c(r_*)$.
 The curve $\Sigma^{(\infty)}$ corresponds to a nonlocal bifurcation and is defined
 below by relation \eqref{eq_SigmaInf}.  }
\label{fig03}
\end{figure*}

In order to illustrate what values $L$ and $E$ take in this case, we define
the plane of first integrals
$$
\mathbb{R}^2_{L,E}=\{ (L, E) \ | \ L>0, E>0 \}.
$$
If  condition \eqref{eq002} is satisfied, then, solving the system of equations  \eqref{eq001} and
\eqref{eq_05}, we obtain the following curve on the plane $\mathbb{R}^2_{L,E}$:
$$
\begin{gathered}
\Sigma^{(r)}= \left\{ (L,E) \ | \ L=L(r_c) , \ E=E(r_c), \ r_c\in(r_1,r_2) \right\}, \\
L(r_c)=\frac{r_c}{b}\sqrt{ \frac{4\big(b - 1\big)\big(r_c - 1\big) +b }{Z(r_c)} }, \quad
E(r_c)=b(r_c-1) \sqrt{\frac{4 - B^2r_c^2}{2r_cZ(r_c)}}.
 \end{gathered}
$$
The curve $\Sigma^{(r)}$ is a bifurcation curve. The values of the integrals on it
correspond to a bifurcation of the function $U(r, \theta)$. As is well known,
using this curve, one can clearly describe the stability of fixed points (for details,
see \cite{BolsinovBifurcation2012}). Indeed, a typical view of the curve $\Sigma^{(r)}$
is presented in Fig. \ref{fig03}, which clearly shows a {\it cusp} type singularity.
Let us introduce the function
\begin{align}
\label{eqNc}
N(r_c)=&64(r_c - 3) + (4Br_c)^2(16r_c^2 - 39r_c  + 21)  
-4(Br_c)^4(24r_c^2 - 51r_c + 25)  +  (Br_c)^6(24r_c^2 - 37 r_c + 15) \nonumber,
\end{align}
then the value $r_c=r_*$, which corresponds to the cusp, satisfies the equation $N(r_*)=0$.
We show that it separates unstable  and stable fixed points. 
It is given by the following matrix:
$$
d^2U=\left(\begin{array}{cc}  \dfrac{N(r_c)}{64r_c^2Z(r_c)b^2(r_c - 1)} & 0 \\0 & \dfrac{4 - B^2r_c^2}{4Z(r_c)} \end{array}\right).
$$
For $r_c\in(r_1, r_*)$ this matrix is negative definite.  Therefore, at these critical points $U$ has a strict maximum.
 Conversely, if $r_c\in(r_*, r_2)$, then the Hessian is positive definite,  and so $U$ has a strict minimum. 
 Also, the first two terms in the Hamiltonian function \eqref{eqRed} are a positive definite quadratic form in momenta. As a result, the maximum point at which $p_r=0$, $p_\theta=0$ corresponds to an unstable fixed point  of the system \eqref{eqRed},  which is of
{\it saddle-center} type, and the minimum corresponds to a stable fixed point of  {\it center-center} type. 

We note that the functions $L(r_c)$ and $E(r_c)$, which parametrize the curve
$\Sigma^{(r)}$, and the inequalities defining the stability of the corresponding solutions
were presented previously in \cite{Lim2015Motion}, but there no visualization and no 
analysis of this curve were presented.

\begin{figure}[!ht]
\includegraphics[scale=1.0]{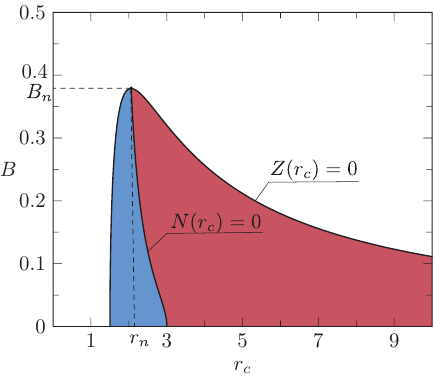}
\caption{On the plane ($r_c, B$), red denotes stable fixed points, and blue indicates, respectively, unstable ones. }
\label{fig02}
\end{figure}

\begin{rem}
Singular points of the curve $\Sigma^{(r)}$ simultaneously satisfy two equations
$$
\frac{dL(r_c)}{dr_c}=0, \quad \frac{dE(r_c)}{dr_c}=0.
$$
In this case they reduce to one equation, $N(r_c)=0$, which defines $r_*$. The proof that
the singular point $r_c = r_*$ is a cusp is fairly standard and so we skip the proof here.
\end{rem}

The values which $r_c$ can take for stable and unstable fixed points depending on the parameter $B$ are presented in Fig. \ref{fig02}. As shown above, their boundaries are defined by the curves $Z(r_c)=0$ and $N(r_c)=0$, from which we obtain the following.

{\it If $B=0$, then $r_*$ takes the largest value $r_*=3$, which corresponds to the innermost stable circular orbit  in the Schwarzschild metric. The magnetic field for $B<B_n$ leads to a decrease in
 $r_*$, but it always holds that $r_*>r_n$, where $r_n=\dfrac{8 + \sqrt{19}}{6}\approx 2.06 $. If $B>B_n$, then there are no relative equilibria.}

\begin{rem}
The critical value obtained above for the parameter $B_n$, after which there are no relative 
equilibria, remains unchanged also for trajectories of light \cite{Galtsov1978Black}.
\end{rem}

It follows from relations \eqref{eqRed} and \eqref{eq_h_level} that the coordinates of the material point always lie in the region
$$
\mathcal{D}^2=\left\{ (r, \theta) \ | \ U(r, \theta) + \dfrac{1}{2} \leqslant 0 \right\},
$$
which we will call, by analogy with the bounded three-body problem, the {\it Hill region}  \cite{ArnoldMathematical2006} or the region of possible motion.
After the values of the first integrals corresponding to the critical points have crossed the curve $\Sigma^{(r)}$), region $\mathcal{D}^2$ transforms, i.e., it undergoes the so-called local bifurcations.
In addition, bifurcations of region $\mathcal{D}^2$ can be due to changes in the behavior of the function $U$ for $r\to + \infty$, which corresponds to a nonlocal bifurcation.

In order to describe them, we introduce the coordinates
$$
y=r \sin\theta, \quad z=r\cos\theta.
$$
Let $(y_c, z_c)$ be the coordinates of the critical points. Then
\begin{equation}
\begin{gathered}
\label{eq_Uyz}
U(y_c, z_c) =\frac{\beta^2 L^2}{2y_c^2} - \left( 1 - \frac{1}{\sqrt{y_c^2 + z_c^2}} \right)^{-1}\frac{E^2}{2 \beta^2}=-\frac{1}{2}, \\
\beta =  1 + \frac{B^2y_c^2}{4}.
\end{gathered}
\end{equation}
In the new coordinates the system \eqref{eqR_Th} can be represented as
\begin{equation}
\begin{gathered}
\label{eq_yz}
\left.\frac{\partial U}{\partial y}\right|_{\substack{y=y_c \\ z=z_c}}=  \frac{\beta(\beta - 2)L^2}{y_c^3} +  (\sqrt{y_c^2 + z_c^2} - 1)^{-2}\frac{ y_c E^2}{2\beta^2\sqrt{y_c^2 + z_c^2}} \\ + \frac{2(\beta -1)\sqrt{y_c^2 + z_c^2}E^2}{\beta^3 y_c(\sqrt{y_c^2 + z_c^2} - 1)}=0  \\
\left.\frac{\partial U}{\partial z}\right|_{\substack{y=y_c \\ z=z_c}}= \left( 1 - \frac{1}{\sqrt{y_c^2 + z_c^2}} \right)^{-2}  \frac{ z_c   E^2}{2\beta^2(y_c^2 + z_c^2)^{\fracs{3}{2}}} = 0.
\end{gathered}
\end{equation}
In the limit $z_c\to \infty$ the second equation is satisfied identically and the left-hand sides of the first equation and of equation \eqref{eq_Uyz} remain constant. Solving them for the values of the first integrals, we obtain
another bifurcation curve on the plane $\mathbb{R}^2_{L,E}$:
\begin{equation}
\label{eq_SigmaInf}
\begin{gathered}
\Sigma^{(\infty)}=\left\{ (L, E)  \  | \  L=L(y_c) \ E=E(y_c), \ y_c\in\left(0, \dfrac{2}{B\sqrt{3}} \right) \right\} \\
L(y_c) = \left( \frac{1}{4} + \frac{1}{B^2y_c^2}\right)^{-1}\frac{\sqrt{2}}{\sqrt{4-3B^2y_c^2}}, \\  E(y_c) = \beta\sqrt{\frac{4 - B^2y_c^2}{4 - 3B^2y_c^2}}.
\end{gathered}
\end{equation}
The curve $\Sigma^{(\infty)}$ also describes bifurcations of regions $\mathcal{D}^2$ for the Melvin solution. The function $L(y_c)$ for the Melvin solution is found in \cite{Melvin1966Orbits}.

\begin{figure*}[!ht]
\includegraphics[scale=0.8]{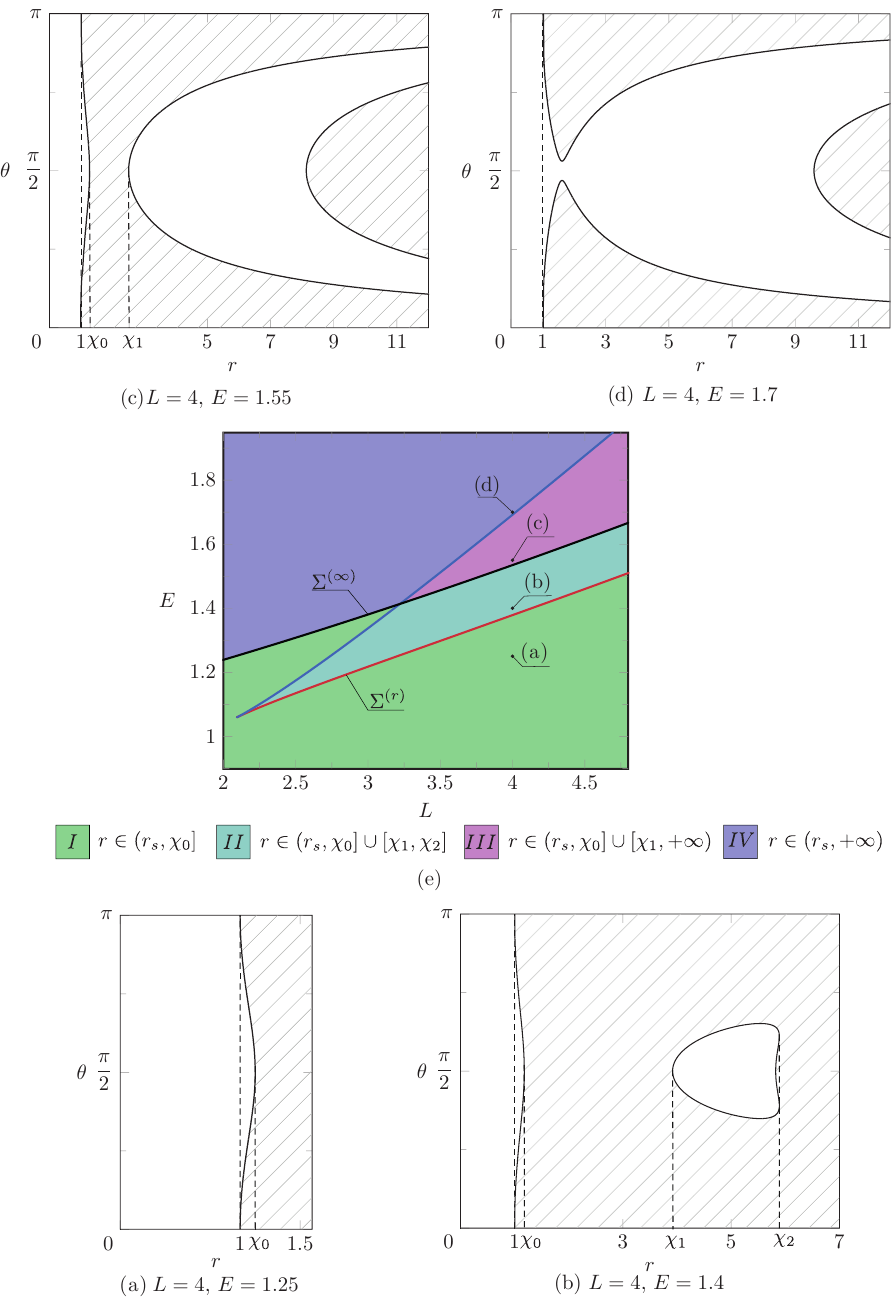}
\caption{ Bifurcation diagram and the Hill regions $\mathcal{D}^2$ shown in white for $r>1$, depending on the values of the first integrals  for the fixed $B=0.15$. Red denotes the part of the curve $\Sigma^{(r)}$ that corresponds to stable fixed points,  and blue indicates the part corresponding to unstable ones. }
\label{fig04}
\end{figure*}

 The relative position of two curves $\Sigma^{(r)}$ and $\Sigma^{(\infty)}$ for different values of the parameter $B$ is shown in Fig. \ref{fig03}. These curves divide the entire plane $\mathbb{R}^2_{L,E}$
into several regions, with different Hill regions $\mathcal{D}^2$.
The number of such regions changes depending on the parameter $B$. Their maximal value is equal to four (see Fig. \ref{fig03}b); in this case a typical view of the Hill regions $\mathcal{D}^2$ is shown in Fig. \ref{fig04}.
Each such Hill region contains a region of possible values of the radial coordinate $r$:
\begin{itemize}
\item[---] I: $r\in(1, \chi_0]$ (see Fig. \ref{fig04}a);
\item[---] II: $r\in(1, \chi_0]\cup[\chi_1, \chi_2]$  (see Fig. \ref{fig04}b);
\item[---] III: $r\in(1, \chi_0]\cup[\chi_1, +\infty)$  (see Fig. \ref{fig04}c);
\item[---] IV: $r\in(1, +\infty)$  (see Fig. \ref{fig04}d).
\end{itemize}
In what follows, as the parameter $B$ increases,  region $III$ disappears 
(see Fig. \ref{fig03}c), and then region $II$ does so when $B>B_n$ (see Fig. \ref{fig03}d).

\section{Plunging trajectories and regularization}
\label{regularization}
In the previous section it was shown that for any values of the integrals $L$ and $E$ the region of possible values of the radial coordinate $r$ includes the event horizon $\mathcal{S}_h$, on which $r=1$.
As can be seen, in this case the system  \eqref{eqRed} has a singularity because $\Gamma \to 0$,
which complicates, for example, the numerical solution of the system \eqref{eqRed} for trajectories reaching the neighborhood of the event horizon $\mathcal{S}_h$.
In order to describe such trajectories, we perform a {\it regularization}
(for details, see \cite{BizyaevInvariant2018, McGeheeTriple1974, Knauf2023}). For this, we determine the momentum $P_r$ and rescale time as
$$
P_r = \Gamma p_r, \quad d\sigma = \frac{d\tau}{\Gamma}.
$$
Then the equations of motion can be represented as
\begin{equation}
\label{eq_reg}
\begin{gathered}
\frac{dP_r}{d\sigma} = -\Gamma^2\frac{\partial \widetilde{H}}{\partial r}, \quad
\frac{dp_\theta}{d\sigma} = -\Gamma\frac{\partial \widetilde{H}}{\partial \theta}, \\
\widetilde{H} = \frac{P_r^2}{2\Lambda^2\Gamma} + \frac{p_\theta^2}{2\Lambda^2 r} + U, \\
\frac{dr}{d\sigma} = \frac{\Gamma}{\Lambda^2}P_r, \quad \frac{d\theta}{d\sigma} = \frac{\Gamma}{\Lambda^2r^2}p_\theta.
\end{gathered}
\end{equation}
This system has the invariant measure
$$
\mu=\frac{1}{\Gamma^2}dP_rdp_\theta dr d\theta.
$$
As can be seen, the density of this measure has a singularity on the invariant manifold  of the system \eqref{eq_reg}
\begin{equation}
\label{eq_I3}
\mathcal{I}^3=\{ (P_r, p_\theta, r, \theta) \ | \ r=1 \}.
\end{equation}
Trajectories on $\mathcal{I}^3$ do not correspond to any real motion of the material point.
However, their study turns out to be useful for the analysis of motion in a neighborhood of the event horizon  $\mathcal{S}_h$.

The equations on the invariant manifold \eqref{eq_I3} have the form
\begin{equation}
\label{eq_Sing}
\begin{gathered}
\frac{dP_r}{d\sigma} = \frac{P_r^2 - E^2}{2\Lambda^2_s}, \quad \frac{dp_\theta}{d\sigma} = \frac{B^2}{2\Lambda_s^3}(P_r^2 - E^2)\sin\theta\cos\theta, \\ \frac{d\theta}{d\sigma}=0, \
 \Lambda_s = 1 + \frac{B^2}{4}\sin^2\theta.
 \end{gathered}
\end{equation}
As can be seen, the first equation decouples from this system. Integrating it explicitly, we find
$$
P_r(\sigma) = E\frac{P_0 - \tanh\left( \dfrac{E\sigma}{2 \Lambda_s^2} \right) }{1 - P_0\tanh\left( \dfrac{E\sigma}{2\Lambda_s^2} \right)}, \quad
P_r(0)=EP_0,
$$
where $P_0$ is the constant defining the initial value of $P_r$. Depending on this constant, the following cases take place:
\begin{itemize}
\item[---] $P_0=1$ and $P_0=-1$ correspond to the fixed points $P_r= E$ and $P_r= -E$; moreover, according to \eqref{eq_Sing}, at these points we have $p_\theta={\rm const}$, $\theta={\rm const}$;
\item[---]  $P_0\in(-1,1)$, then $\lim\limits_{\sigma\to + \infty} P_r =  - E$ and $\lim\limits_{\sigma\to - \infty} P_r =  E$;
\item[---]  $P_0>1$, then  $\lim\limits_{\sigma\to  -\infty} P_r =  E$ and $P_r \to +\infty$ in finite time $\sigma = \dfrac{2\Lambda_s^2}{E}\arctanh \dfrac{1}{P_0}$;
\item[---]  $P_0<-1$, then  $\lim\limits_{\sigma\to  +\infty} P_r =  -E$ and $P_r \to -\infty$ in finite time $\sigma = -\dfrac{2 \Lambda_s^2}{E}\arctanh \dfrac{1}{P_0}$.
\end{itemize}

\begin{figure}[!ht]
\includegraphics[scale=1.1]{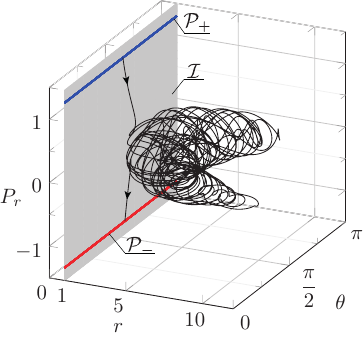}
\caption{An example of the plunging trajectory that tends to $\mathcal{P}_\pm$ as $\sigma \to \mp \infty $ for fixed  $B=0.15$, $L=2.8$, $E=1.3$. For the initial conditions $P_r(0)=1.2842$, $p_\theta(0)=0.54629$, $r(0)=1.00441$, $\theta(0)=1.63352$.}
\label{fig05}
\end{figure}

We now consider the regularized system \eqref{eq_reg}. It possesses two families of fixed points
$$
\mathcal{P}_\pm^2 =\{ P_r=\pm E, \ p_\theta=p_*={\rm const}, \ r=1, \ \theta=\theta_*={\rm const}  \}.
$$
Throughout the remainder of this section, the upper sign corresponds to $\mathcal{P}^2_+$, and the lower sign, to $\mathcal{P}^2_-$.
As shown above, the trajectories lying on $\mathcal{I}^3$ approach the families $\mathcal{P}_\pm^2$  in time $\sigma\to \mp\infty$. We show that this also holds true for trajectories  in a neighborhood of $\mathcal{I}^3$.
To do so, we represent the linear system in a neighborhood of these equilibrium points in the following form:
\begin{equation}
\label{eq_Linear}
\begin{gathered}
\frac{dz_1}{d\sigma} = \lambda z_1, \quad \frac{dz_3}{d\sigma} = \lambda z_3, \quad  \frac{dz_4}{d\sigma} = \frac{p_*}{\Lambda_*^2}z_{3},\\
\frac{dz_2}{d\sigma} = \sin ({2\theta_*}) \frac{\lambda B^2 }{2\Lambda_*} z_1 + \cos \theta_*\left[  \frac{2p_*^2{\sin \theta_*}(\Lambda_* - 1)}{\Lambda^3_{*}} -  \Lambda_*(\Lambda_*-2)L^2 \right] \frac{z_3}{\sin^3\theta_{*}}, \\
 \Lambda_* = 1 + \frac{B^2}{4}\sin^2\theta_*,
\end{gathered}
\end{equation}
where $\boldsymbol{z}=(P_r\mp E, p_\theta - p_*, r-1, \theta - \theta_*)$ is the deviation from the fixed points and the following parameter has been introduced:
\begin{align*}
\mathcal{P}_+^2:  \ \lambda &= \frac{E}{\Lambda_*^2}>0, \\
\mathcal{P}_-^2:  \ \lambda &= -\frac{E}{\Lambda_*^2}<0.
\end{align*}
The first two equations decouple from the system \eqref{eq_Linear}. Integrating them explicitly, we find
\begin{equation}
\label{eq_z12}
z_1 = C_1 e^{\lambda \sigma}, \quad z_3 = C_3 e^{\lambda \sigma}.
\end{equation}
Thus, for trajectories in a neighborhood of $\mathcal{P}^2_+$ we have $r \to 1$ as $\sigma \to -\infty$, and vice versa, in a neighborhood of $\mathcal{P}^2_-$ we have $r \to 1$ as $\sigma \to +\infty$. We show that the initial proper time of the particle $\tau$ takes a finite value in this case. Indeed, using the solution of the linear system \eqref{eq_z12}, we find
$$
\tau = \int\limits_0^{\sigma} \left(1 - \frac{1}{1 + C_3 e^{\lambda \tilde{\sigma}}}\right) d\tilde{\sigma} = \frac{1}{\lambda} \ln \frac{1 + C_3e^{\lambda \sigma}}{1 + C_3}.
$$
For $\mathcal{P}_\pm^2$, in the case $\sigma \to \mp\infty$ this time $\tau$ takes a finite value.
Next, we will call the trajectories of the system \eqref{eq_reg} for which $r\to 1$ the 
{\it plunging trajectories}.
An example of such a trajectory is presented in Fig. \ref{fig05}.
\begin{rem}
We rewrite the linear system \eqref{eq_Linear} as $\dfrac{d\boldsymbol{z}}{d\sigma}={\bf A}\boldsymbol{z}$, where the matrix {\bf A} has two zero eigenvalues. As a result, there exist vectors $\boldsymbol{b}$ satisfying  the equation ${\bf A}^{T}\boldsymbol{b}=0$. It can be shown that each such vector corresponds to the linear integral $F=(\boldsymbol{b},\boldsymbol{z})$.  In this case there exist two linearly independent integrals  
$$
\begin{gathered}
F_1 = (\boldsymbol{b}_1, \boldsymbol{z}), \quad F_2 = (\boldsymbol{b}_2, \boldsymbol{z}), \\
\boldsymbol{b}_1 = \left(0, 0, \mp p_*, E \right), \\
\boldsymbol{b}_2 = \left( - \frac{B^2E}{2\Lambda_*}\sin 2\theta_*, E,  \pm \frac{\cos \theta_*}{ \sin\theta_*}\left( \frac{\Lambda_*^3(\Lambda_* - 2)}{\sin^2\theta_*}L^2 \right) , 0 \right).
\end{gathered}
$$
Fixing the values of these integrals and using the solution  \eqref{eq_z12}, we can easily obtain $z_2(\sigma)$ and  $z_4(\sigma)$.
\end{rem}

\section{Bounded trajectories and a Poincar\'{e} map}
\label{sectionMap}
The system \eqref{eqRed} possesses the discrete symmetry
\begin{equation}
\label{symmetry}
p_\theta \to -p_\theta, \quad \theta \to \pi - \theta.
\end{equation}
The fixed points of this symmetry define the invariant submanifold of the system \eqref{eqRed}:
$$
p_\theta =0, \ \theta=\dfrac{\pi}{2},
$$
on which plane trajectories lie. They are analyzed in detail in  \cite{Galtsov1978Black, Dadhich1979Trajectories, Galtsov1986Book, Esteban1984}. Next, we consider the spatial trajectories.
To visualize them, we will construct a Poincar\'{e} map for the regularized system \eqref{eq_reg}. We describe its construction.
\begin{itemize}
\item[---]  As a secant of the vector field of the system \eqref{eq_reg} we choose
a manifold given by the relation
\begin{equation}
\label{SectionMap}
\theta=\dfrac{\pi}{2},
\end{equation}
which corresponds to the equatorial plane.
\item[---] We will construct the Poincar\'{e} map in the variables $(r, P_r)$.
Therefore, from the given initial values of $r$ and $P_r$ lying on the secant, we will determine the remaining initial value for the momentum $p_\theta$ from the Hamiltonian function $\widetilde{H}=-\dfrac{1}{2}$
as follows:
$$
\begin{gathered}
p_\theta=\Lambda\sqrt{2r\Delta(r, P_r)}, \\ \Delta(r, P_r) =  \frac{8(E^2 - P_r^2)}{\Gamma(4 + B^2r^2)^2} - \frac{1}{2} - \frac{L^2}{32r^2}(4 + B^2r^2)^2.
\end{gathered}
$$
As can be seen, for the initial $r$  and $P_r$ the condition $\Delta(r, P_r)>0$ must be satisfied.
\item[---]
For the above initial conditions we perform numerical integration of the system \eqref{eq_reg}.
Simultaneously, using the Henon method \cite{HenonNumerical1982}, we will find intersections of the
trajectories with the given section, for which $\dot{\theta}>0$.
Depicting the found values of the coordinates on the plane $(r, P_r)$, we obtain a two-dimensional
Poincar\'{e} point map.
\end{itemize}

\begin{rem}
As is well known, the fixed points of the Poincar\'{e} map correspond to the periodic solutions of the initial system, and the invariant curves, to the quasiperiodic solutions lying on a torus.
\end{rem}

We note that the plunging trajectories, $r \to 1$, or the trajectories going to infinity, $r\to +\infty$, can cease to intersect the chosen secant. Therefore, the Poincar\'{e} map can be used only in the case where the recurrence of the trajectories is observed.
In what follows, we will call the trajectories for which the radial coordinate always remains bounded and does not tend to $1$ {\it bounded trajectories}. They exist for the values of the integrals $L$ and $E$ from region $II$  (see Fig. \ref{fig04}) or possibly for other values of the integrals, but near stable periodic solutions.

\begin{figure*}[!ht]
\includegraphics[scale=0.95]{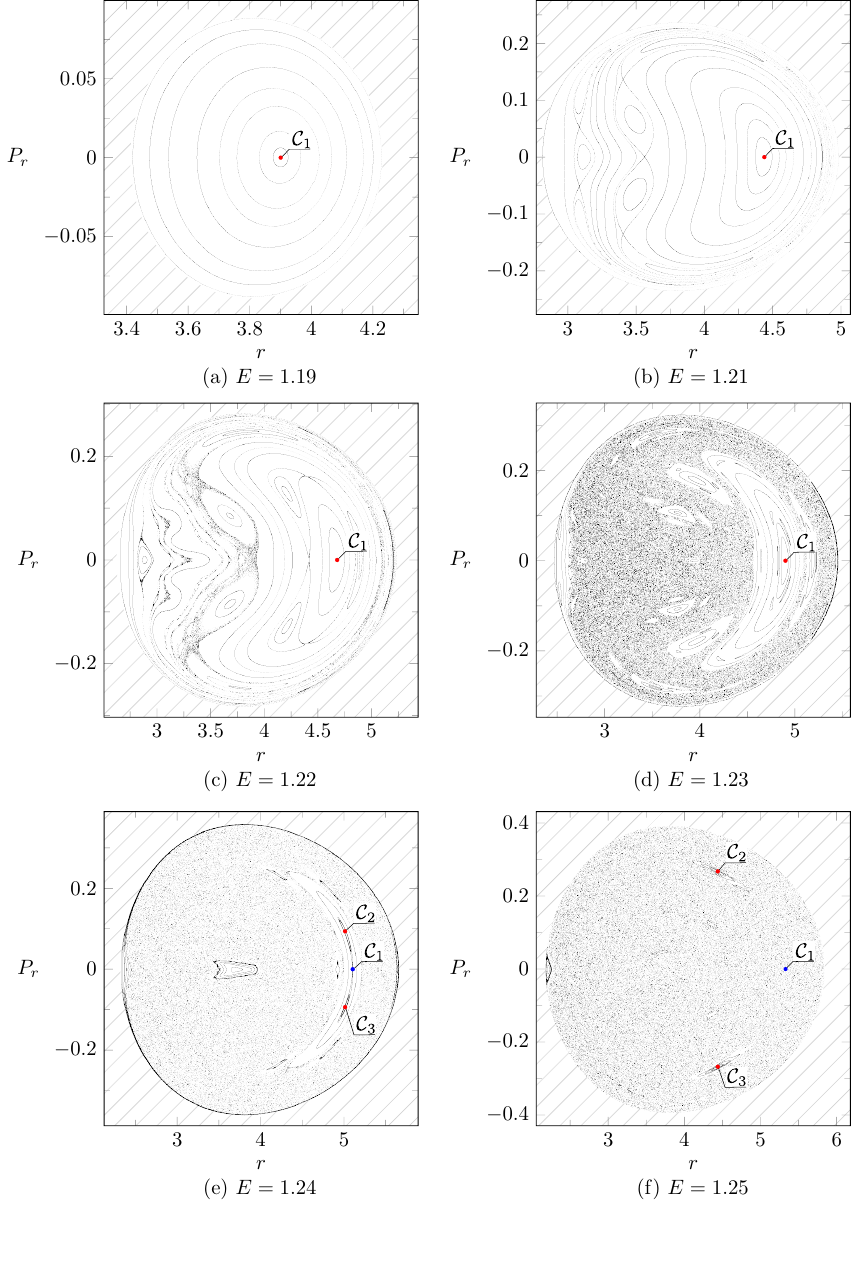}
\caption{A Poincar\'{e} map for fixed $B=0.15$, $L=2.8$  and different values of
$E$. Red denotes the stable fixed points, and blue indicates the unstable ones.
The hatched area is the region in which  $\Delta(r, P_r)<0$, i.e., which does not
belong to the domain of definition of the Poincar\'{e} map. }
\label{figGrup1}
\end{figure*}

We first consider the first case.
For the values of the integrals $L$ and $E$, which correspond to the stable part of the curve $\Sigma^{(r)}$ (i.e., to the lower boundary of region $II$) the system \eqref{eqRed} possesses a fixed point of center-center type. Above this part of the curve $\Sigma^{(r)}$, a periodic solution appears in the system. In the Poincar\'{e} map, this solution corresponds to the fixed point $\mathcal{C}_1$ (see Fig. \ref{figGrup1}a), which in this case is surrounded by invariant curves
 (because these trajectories lie in a neighborhood of a fixed point  of center-center type).  Next, as energy $E$ increases, the invariant curves break down 
(see Fig.~\ref{figGrup1}b,c,d). Then the
stable fixed point undergoes a pitchfork bifurcation. After that, it becomes unstable and in its neighborhood two new stable fixed points  $\mathcal{C}_2$ and $\mathcal{C}_3$ appear (Fig. \ref{figGrup1}e). 
Then these stable points move away from each other and lose stability
after a period-doubling bifurcation.

\begin{figure*}[!ht]
\includegraphics[scale=0.65]{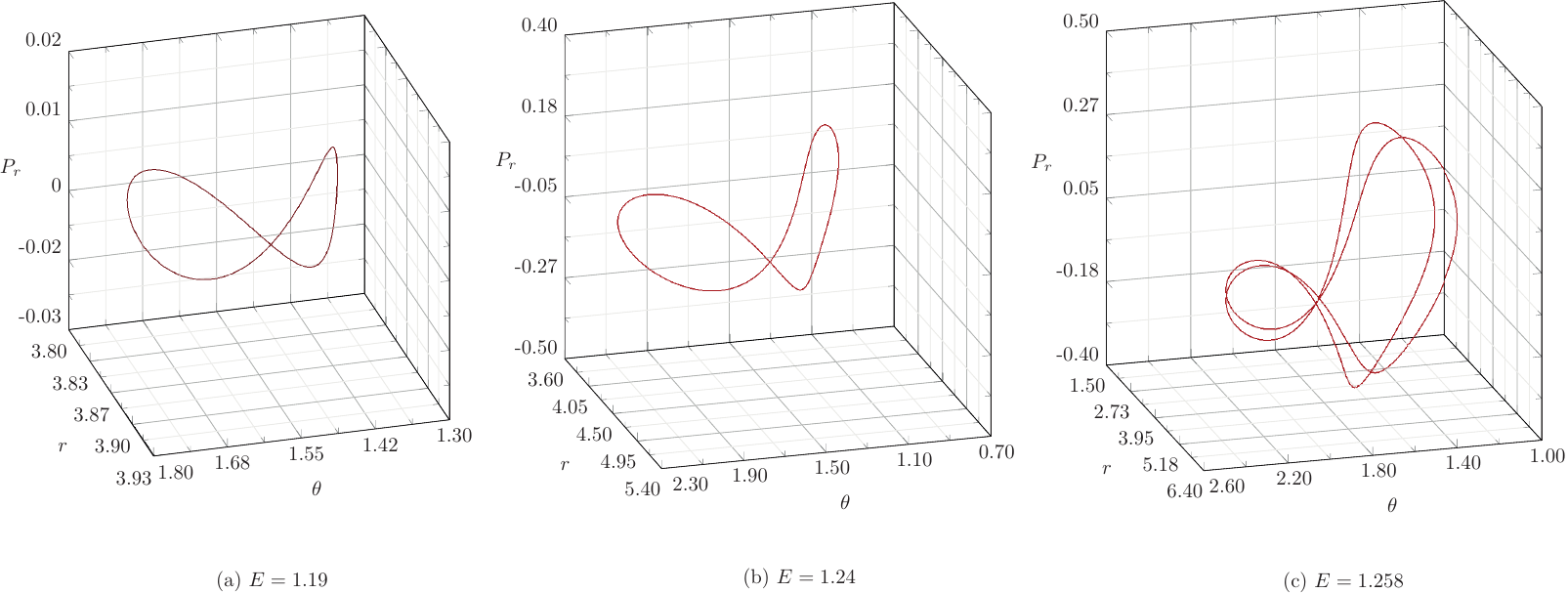}
\caption{ Projection of the stable periodic solution of the regularized system
\eqref{eq_reg}  for fixed parameters $B=0.15$, $L=1.28$ and different values of
energy $E$. The initial conditions are given as follows: (a) $P_r(0)=0$, $p_\theta(0)=0.3932$, $r(0)=3.9005$, $\theta(0)=\dfrac{\pi}{2}$, (b) $P_r(0)=0.0938$, $p_\theta(0)=1.3962$, $r(0)=5.0135$, $\theta(0)=\dfrac{\pi}{2}$, (c) $P_r(0)=0.3331$, $p_\theta(0)=1.1167$, $r(0)=3.8058$, $\theta(0)=\dfrac{\pi}{2}$.  }
\label{figGrup3}
\end{figure*}

\begin{figure}[!ht]
\includegraphics[scale=0.9]{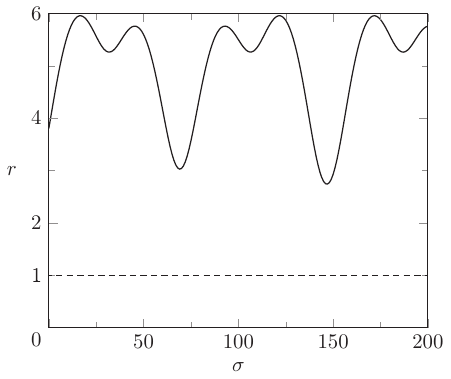}
\caption{ Dependence $r(\sigma)$ for initial conditions and parameters corresponding to
Fig.~\ref{figGrup3}c.  }
\label{Curve_r_sigma}
\end{figure}

A few examples of trajectories of the system that correspond to the stable fixed point
of the Poincar\'{e} map are presented in Fig. \ref{figGrup3}. For example, Fig.
\ref{figGrup3}a shows a trajectory corresponding to the fixed point shown in red in Fig.
\ref{figGrup1}a, this periodic solution is symmetric about the transformation
\eqref{symmetry}. Next, the pitchfork bifurcation described above gives rise to an
asymmetric periodic solution (Fig. \ref{figGrup3}b). Figure \ref{figGrup3}c shows an example of
a periodic solution that intersects the Poincar\'{e} section twice (i.e., the fixed point
of period 2).  This trajectory can  get close to the event horizon $\mathcal{S}_h$
(see Fig. \ref{Curve_r_sigma}).


\begin{figure}[!ht]
\includegraphics[scale=0.9]{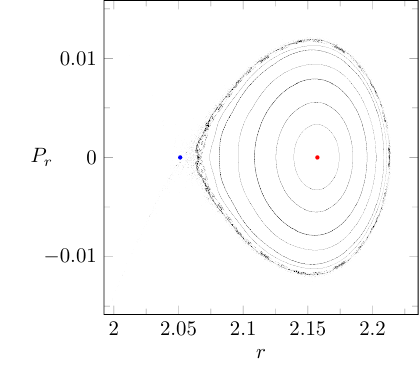}
\caption{A Poincar\'{e} map for fixed $B=0.4>B_n$, $L=2.8$ and $E=1.815$.
Red denotes the stable fixed point, and blue indicates the unstable one.}
\label{periodicSolution}
\end{figure}

\begin{figure}[!ht]
\includegraphics[scale=0.85]{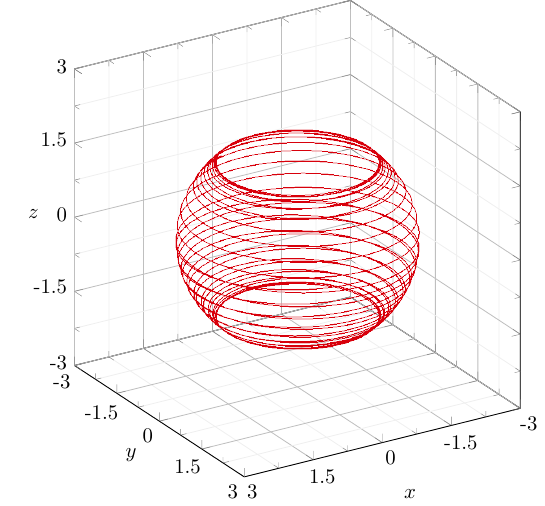}
\caption{Trajectory in the space $x=r \sin \theta\cos\varphi$, $y=r\sin \theta \sin \varphi$,
$z=r \cos \theta$ for the fixed $B=0.4>B_n$, $L=2.8$, $E=1.815$ of the initial conditions
in the neighborhood of the stable fixed point in Fig. \ref{periodicSolution}: $P_r(0)=0$,
$r(0)=2.15737$, $\theta(0)=\dfrac{\pi}{2}$, $p_\theta(0)=2.55158$.  }
\label{Curve3d}
\end{figure}

It was shown in Section \ref{HillReg} that for a sufficiently strong magnetic field
$B>B_n$ there are no relative equilibria. This does not imply the absence of
other bounded trajectories. Such trajectories can arise after a saddle-node bifurcation.
As a result, a pair of fixed points (a stable and an unstable one) arises in the map.
An example of the Poincar\'{e} map after such a bifurcation is presented in Fig.
\ref{periodicSolution}. Figure \ref{Curve3d} shows a trajectory with initial
conditions in the neighborhood of the stable fixed point.

\section{Separatrix splitting}
\label{separatrix}
As noted above, if $B=0$ then we obtain the Schwarzschild case. In this case,
the equations of motion \eqref{eqRed} possess the additional integral
\begin{equation}
\label{eqFIntegral}
F=p_\theta^2 + \frac{L^2}{\sin^2\theta},
\end{equation}
and hence are integrable by the Liouville-- Arnold theorem.
A detailed analysis of the trajectory of the material particle is performed, for example,
in \cite{HagiharaTheory1930, Chandrasekhar1998Mathematical}.
Let us consider in more detail the periodic solution which, after rescaling time as
$dt=r^2du$, can be represented as
\begin{equation}
\label{eq_sw_per}
\begin{gathered}
r=r_u={\rm const}, \quad p_r=0, \\
\cos\theta(u) = \sqrt{1 - L^2 a^2} \cos\left( \frac{u}{a} +  C_\theta   \right) , \\
p_\theta(u) = \frac{\sqrt{2}}{a} \sin\left( \frac{u}{a} +  C_\theta   \right)\left(\frac{1 + L^2a^2}{1 - L^2a^2} - \cos \left(\frac{2u}{a} +  2C_\theta\right)    \right)^{-\fracs{1}{2}},
\end{gathered}
\end{equation}
where the constant $C_\theta$ defines the initial position of angle $\theta$ and it is
assumed that for the projection of the momentum the condition  $L<\dfrac{1}{a}$ is
satisfied. In addition, energy $E$ and parameter $a$ depend on $r_u$ as follows:
$$
E=\frac{\sqrt{2}(r_u - 1)}{\sqrt{r_u(2r_u - 3)}}, \quad a=\frac{\sqrt{2r_u - 3 }}{r_u}.
$$

If $r_u\in(2, 3)$, then the solution \eqref{eq_sw_per} is unstable and there exist
doubly asymptotic trajectories (separatrices) that approach it as $u\to~\pm\infty$, i.e.,
we have homoclinic trajectories. The solution for these trajectories can be obtained in
elementary functions
\begin{equation}
\label{eq_solution}
r(u) = \frac{r_u \cosh^2 (\gamma u) }{\dfrac{r_u}{r_a} + \sinh^2(\gamma u)}, \quad
p_r(u) = - \frac{2r_u\tanh (\gamma u) }{1 + \dfrac{r_u(r_u - 1)}{r_a - r_u}},
\end{equation}
where the following parameters have been introduced:
$$
r_a = \frac{r_u}{r_u-2}>r_u, \quad \gamma = \frac{\sqrt{r_u( 3 - r_u)}}{2\sqrt{2 r_u - 3 }}.
$$
The functions $\theta(u)$ and $p_\theta(u)$ for homoclinic trajectories are also defined
by relations \eqref{eq_sw_per}. They differ in the value of the parameter
$C_\theta$. The solution in the form \eqref{eq_solution} can be obtained using
\cite{BizyaevBifurcation2022}, where doubly asymptotic plane trajectories
for the Kerr metric are obtained.

In addition to the homoclinic trajectories described above there exist two other types of
trajectories asymptotic to the unstable periodic solution \eqref{eq_sw_per}. We denote them
by $\mathcal{S}_s$ and $\mathcal{S}_u$.
The trajectories corresponding to $\mathcal{S}_u$ approach the solution \eqref{eq_sw_per}
as $u\to-\infty$, and reach the event horizon in finite time $u$, while the trajectories
$\mathcal{S}_s$  ``start'' on the horizon and asymptotically approach the periodic solution
\eqref{eq_sw_per} as $u\to+\infty$.

\begin{rem}
 If $r_u<2$, then $r_a<0$ and the solution  \eqref{eq_solution}
does not describe a homoclinic trajectory since in this case  $r\to +\infty$ (i.e., the
separatrices go to infinity).
 If $r_u>3$, then the solution \eqref{eq_sw_per} is stable.
\end{rem}

Consider in more detail the homoclinic trajectories described above. These trajectories are asymptotic manifolds (stable and unstable manifolds) which coincide in the case $B=0$ considered above.  After perturbation  (i.e.,  $B\neq0$)  these manifolds persist, but can intersect. They intersect if the Melnikov integral has isolated zeros  \cite{WigginsBook, KozlovIntegrability1983}. 

We represent the Melnikov integral in the following form (see, e.g., \cite[p. 272]{RobinsonDynamical1995}, \cite{Robinson1988}):
$$
\begin{gathered}
J(C_\theta) = \int\limits_{-\infty}^{{+\infty}} \frac{dF_0}{du}du, \\
 F_0 = F - 2r_u^2 \left.H \right|_{B=0}=\frac{r r_u^2E^2}{r - 1} + \frac{r^2 - r_u^2}{r^2}\left( \frac{L^2}{\sin^2\theta} + p_\theta^2  \right) - \frac{r - 1}{r}r_u^2p_r^2,
\end{gathered}
$$
where $F$ is the integral in the Schwarzschild case \eqref{eqFIntegral}.
As can be seen,  $F_0$ is also an additional integral  of equations  \eqref{eqRed}, but in contrast to $F$, for $F_0$,
  at the points 
of the unstable periodic solution \eqref{eq_sw_per} we have $dF_0=0$, which is necessary 
for convergence of the Melnikov integral \cite{KozlovIntegrability1983,Robinson1988}.
In addition, in this integral the differentiation of $F_0$ has been performed along the 
trajectory of the system \eqref{eqRed} using the new time variable, and integration 
has been performed along the homoclinic solutions described above. 
Taking this into account, we represent the additional integral as
$$
\begin{gathered}
J(C_\theta) = B^2J^{(1)}(C_\theta) + O(B^4), \\
J^{(1)}(C_\theta) = r_u^2\int\limits_{-\infty}^{+\infty} r^2(u)\sin^2 \theta(u) \big(r(u) - 1\big) p_r(u)du  - \frac{1}{2}\int\limits_{-\infty}^{+\infty}r^2(u)\sin\big(2\theta(u)\big)\big( r^2(u) - r^2_u  \big)p_\theta(u)du.
\end{gathered}
$$
In this case the integral for $J^{(1)}(C_\theta)$ can be expressed in elementary functions 
and has the following form:
$$
\begin{gathered}
J^{(1)}(C_\theta) = -\frac{\pi r_u^5 (1 - a^2L^2) f(r_u)}{4(r_u - 2)^3 \sinh\dfrac{\pi}{a\gamma}}  \sin (2C_\theta) , \\
f =  5 (2 r_u - 3)^2\cosh\left( \frac{2}{a\gamma} \arccos\dfrac{\sqrt{r_u}}{\sqrt{r_a}}  \right) +
 \frac{8r_u^2(10r_u - 39)  + 358r_u - 103}{4 \sqrt{r_u(r_u - 2)}} \sinh\left( \frac{2}{a\gamma} \arccos\frac{\sqrt{r_u}}{\sqrt{r_a}}  \right).
 \end{gathered}
$$
\begin{figure*}[!ht]
\includegraphics[scale=0.85]{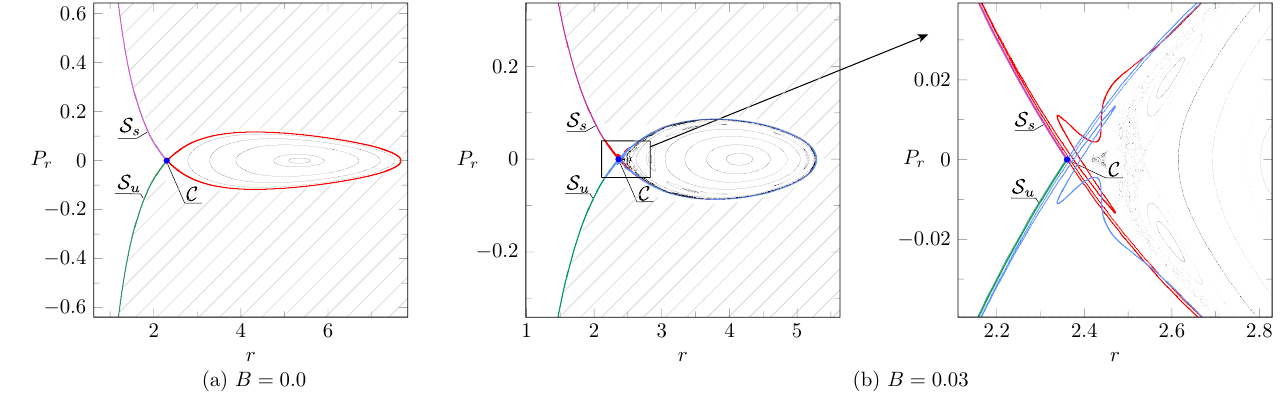}
\caption{A Poincar\'{e} map for fixed parameters $L=1.8$, $E=0.9584$.
The hatched area is the region in which $\Delta(r, P_r)<0$, i.e., which does not belong to 
the domain of definition of the Poincar\'{e} map.
The curves shown in red, blue, green and purple correspond to separatrices.   }
\label{figGrup4}
\end{figure*}

\begin{rem}
In order to obtain an expression for $J^{(1)}(C_\theta)$,
the following relations can be used:
$$
\small{
\begin{gathered}
\int\limits_{-\infty}^{+\infty}\sin(\alpha x + C)\frac{\cosh^4 x (\beta^2 + \cosh(2x))}{(\beta^2 + \sinh^2 x)^4}dx =
 -\frac{\pi\sin(2C)}{16\beta^6 \sinh\dfrac{\alpha\pi}{2}}   \left[ \frac{p_1\sinh( \alpha \arccos\beta) }{\beta\sqrt{1 - \beta^2}}
+ \frac{\alpha}{3} p_2\cosh( \alpha \arccos\beta)\right], \\
p_1(\alpha,\beta)=2\beta^2(\beta^2 - 1)(2\beta^2+1)\alpha^2 - 6\beta^2 - 5, \\ p_2(\alpha,\beta) = \beta^2(\beta^2-1)\alpha^2 - 4 \beta^2(5\beta^2 + 7) - 15,  \\
\int\limits_{-\infty}^{+\infty}\cos(\alpha x + C)\frac{\cosh^3 x \sinh x }{(\beta^2 + \sinh^2 x)^3}dx = 
 -\frac{\pi \alpha\sin(2C)}{8\beta^2 \sinh\dfrac{\alpha\pi}{2}} \left[ \frac{1 + 2\beta^2}{\beta\sqrt{1 - \beta^2}} \sinh( \alpha \arccos\beta)
+\alpha \cosh( \alpha \arccos\beta)\right],\\ \beta<1, \quad \alpha>0.
\end{gathered}
}
$$
They can be obtained using deductions.
\end{rem}

As can be seen, the function $J^{(1)}(C_\theta)$ always has an isolated zero, 
and hence the separatrices do  intersect for $B\neq 0$.
Firstly, this implies that in the general case the system has no additional integral
(for details, see the review \cite{Gelfreich2001Splitting}).
Secondly, there may exist stable periodic trajectories visiting the separatrix holes.
These trajectories are not generated by any periodic trajectories of the unperturbed 
system (i.e., when $B=0$), but it is fairly difficult to numerically find them 
(for details, see  \cite{Simo2011Some,Piftankin2007Separatrix}).

An example of the Poincar\'{e} map  in which the separatrices are numerically constructed 
is presented in Fig. \ref{figGrup4}.  Blue denotes the fixed point $\mathcal{C}$, which 
corresponds to the unstable periodic solution \eqref{eq_sw_per}, green and purple indicate 
the separatrices corresponding to $\mathcal{S}_s$ and $\mathcal{S}_u$. The red curve in 
Fig. \ref{figGrup4}a is a curve obtained by iterations of the homoclinic trajectories, 
and the red and blue lines in Fig.  \ref{figGrup4}b are intersecting stable and
unstable separatrices
to the fixed point $\mathcal{C}$. Plunging trajectories lie in the neighborhood of these 
intersecting separatrices  because, after reaching a neighborhood of point
$\mathcal{C}$m, they approach the event horizon along $\mathcal{S}_u$ as 
$\sigma \to +\infty$, and as $\sigma \to -\infty$ they approach the horizon along 
$\mathcal{S}_s$.  An example of such a trajectory was
presented in Fig. \ref{fig05}.

\section{Conclusion}
This paper is concerned with the trajectories of neutral particles in the Schwarzschild-Melvin space. Using the methods of qualitative analysis from dynamical systems theory, a classification of these trajectories is made according to the values of the first integrals $L$ and $E$ and the parameter characterizing the magnetic field $B$. 

The circular  trajectories $r=r_c$ are considered in detail. For them, the functions $L(r_c)$ and $E(r_c)$ are obtained.
Earlier,  similar relations  and stability conditions were obtained  in \cite{Lim2015Motion, Esteban1984}, but they were not completely analyzed. 
In this paper, it is shown that these functions,   $L(r_c)$ and $E(r_c)$, define the parametrically given curve $\Sigma^{(r)}$ on the plane of first integrals, which has a cusp singularity for $r_c=r_*$. The value $r_c=r_*$ separates stable and unstable circular orbits. In addition, this curve, together with $\Sigma^{(\infty)}$  divides the plane of first integrals into four different regions of possible motion. 
 
If $B>B_n$, where $B_n\approx0.379$,  then there are no circular orbits. Earlier, this result was obtained in \cite{Esteban1984} for neutral particles  and in  \cite{Galtsov1978Black} for light beam trajectories.
In addition, in  \cite{Karas1992Chaotic,Li2019Chaotic} it was shown that, if $B<B_n$, then chaotic trajectories occur. This was shown using a numerically constructed Poincare map and by calculating the largest Lyapunov exponent. In this paper a detailed analysis is made of the bifurcations of periodic solutions using a Poincar\'{e} map. In particular, it is shown that bounded (including chaotic) trajectories can arise after a saddle-node bifurcation for $B>B_n$ as well.  Also, the Melnikov integral is calculated. It shows that, for a small value of the parameter $B$,
stable and unstable manifolds intersect. This is also illustrated by the Poincar\'{e} map.

\begin{acknowledgments}
This work was carried out within the framework of the state assignment of the Ministry 
of Science and Higher Education of Russia (FEWS-2024-0007).
\end{acknowledgments}

\bibliography{Bizyaev}

\end{document}